\documentclass[11pt]{amsproc}

\newtheorem{Theorem}{Theorem}

\newtheorem{theorem}{Theorem}[section]
\newtheorem{lemma}[theorem]{Lemma}
\newtheorem{proposition}[theorem]{Proposition}

\theoremstyle{definition}

\newtheorem{example}[theorem]{Example}

\theoremstyle{remark}

\numberwithin{equation}{section}

\newcommand{\abs}[1]{\lvert#1\rvert}

\newcommand{\R}{{\mathbb R}}

\newcommand{\Z}{{\mathbb Z}}

\newcommand{\Int}{{{\rm Int}\,}}

\newcommand{\RR}{{\mathcal R}}
\newcommand{\EE}{{\mathcal E}}

\begin{document}

\title{Flows of flowable Reeb homeomorphisms}


\author{Shigenori Matsumoto}
\address{Department of Mathematics, College of
Science and Technology, Nihon University, 1-8-14 Kanda, Surugadai,
Chiyoda-ku, Tokyo, 101-8308 Japan
}
\email{matsumo@math.cst.nihon-u.ac.jp
}
\thanks{The author is partially supported by Grant-in-Aid for
Scientific Research (C) No.\ 20540096.}
\subjclass{37E30}

\keywords{Reeb foliations, homeomorphisms, topological conjugacy.}

\date{\today }

\begin{abstract}
We consider a fixed point free homeomorphsim $h$ of 
the closed band $B=\R\times[0,1]$ which leaves each leaf
of a Reeb foliation on $B$ invariant. Assuming $h$ is
the time one of various topological flows, we compare the
restriction of the flows on the boundary.
\end{abstract}

\maketitle

\section{Introduction}

Orientation preserving and fixed point free homeomorphisms 
of the plane are called {\em Brouwer homeomorphisms}.
Since the seminal work of L. E. Brouwer nearly 100 years ago,
they draw attentions of many mathematicians (\cite{K},\cite{HT},
\cite{F},\cite{Fr}
\cite{G}, \cite{N}). Nowadays there still remains interesting problems
about them.

Besides those Brouwer homeomorphisms which are topologically
conjugate to the translation, the  simplest ones are perhaps
those which preserve the leaves of Reeb foliations; the main
theme of the present notes. It is simpler and loses
nothing to consider their restriction to the Reeb component.

Let $B=\R\times[0,1]$ be a closed band, and denote 
$\partial_iB=\R\times\{i\}$ ($i=0,1$) and $\Int B=\R\times(0,1)$.
An oriented foliation $\RR$ on  $B$ is
called a {\em Reeb foliation} if $\partial_0B$ with the positive
orientation and $\partial_1B$ with the negative orientation are
leaves of $\RR$ and the foliation restricted to the interior
$\Int B$ is a bundle foliation. The leaf space of a Reeb foliation 
is homeomorphic to the non Hausdorff space obtained by glueing two 
copies of $[0,\infty)$ along $(0,\infty)$. This shows (\cite{HR})
that any two Reeb foliations are mutually topologically equivalent.

A homeomorphism $h:B\to B$ is called a {\em Reeb homeomorphism} if
$h$ preserves each leaf of a Reeb foliation and $h(x)>x$ for any $x\in
B$,
where $>$ is the total order on a leaf given by the orientation.
In this paper we consider {\em flowable} Reeb homeomorpshisms $h$,
i.\ e.\ those which are the time one of topological flows.
(F. B\'eguin and F. Le Roux constructed in \cite{BL}
 examples of non flowable Reeb homeomorphisms.)
When $h$ is flowable, $h$ can be embedded as the time one
into a great variety of flows. The purpose of this paper is to compare
the restriction of one flow to the boundary of $B$ with that
of another. This problem is motivated by a result of \cite{L2}
which states that if two flows have a common orbit foliation,
then their restrictions to the boundary are the same.

Let
$$P=\{(x,y)\mid x\geq0, y\geq0\}-\{(0,0)\}.$$
Notice that $P$ is homeomorphic to $B$. Consider a homeomorphsim
$h_P:P\to P$ defined by
$$h_P(x,y)=(2x,2^{-1}y)).$$
A homeomorphism of $B$ is called a {\em standard Reeb homeomorphism}
if it is topologically conjugate to $h_P$, and 
{\em nonstandard} otherwise. 
It is known (\cite{BL}) that there are nonstandard flowable 
Reeb homeomorphisms.
The main result of this paper is
the following.

\begin{Theorem} \label{main}
(1) Assume $h$ is a standard Reeb homeomorphsim of $B$. For $i=0,1$,
let $\{\psi_i^t\}$ be an arbitrary flow on $\partial_iB$ 
whose time one is the restriction of $h$.
Then there is a flow $\{\varphi^t\}$ on $B$, an extension of 
both $\{\psi_0^t\}$ and $\{\psi_1^t\}$, whose time one is $h$.

(2) If $h$ is a nonstandard flowable Reeb homeomorphism, 
there is a homeomorphism from $\partial_0B$ to $\partial_1B$ which
 commutes
with any flow whose time one is $h$.

\end{Theorem}

Notice that (1) is immediate from the model $h_P$ which is a
product map.
After we prepare some necessary prerequisites in Sect.\ 2, we prove
Theorem \ref{main} (2) in Sect.\ 3. Sect.\ 4 is devoted to two examples
of nonstandard flowable Reeb homeomorphisms, one for which Theorem
\ref{main} (2) is the optimal, and the other for which the restriction
of the flow to the boundary is unique.

\section{Preliminaries}

Let $\RR$ be a Reeb foliation on $B$. A topological flow on $B$ is called an
$\RR$-{\em flow} if its oriented orbit foliation is $\RR$. 
Let $\EE$ be the set of the topological conjugacy classes of $\RR$-flows.
We shall summerize a main result of \cite{L},
a classification of $\EE$, which will play
a crucial role in what follows.

For $i=0,1$ let
$\gamma_i:[0,\infty)\to B$ be a continuous path such that 
$\gamma_i(0)\in \partial_iB$
and that $\gamma_i$ intersects every interior leaf of $\RR$
at exactly one point. Let us  parametrize
$\gamma_i$ so that for any $y>0$ the points $\gamma_0(y)$
and $\gamma_1(y)$ lie on the same leaf of $\RR$.
Let $\{\Phi^t\}$ be an $\RR$-flow.
Then one can define a continuous function 
$$f_{\{\Phi^t\},\gamma_0,\gamma_1}:(0,\infty)\to \R$$
by setting that $f_{\{\Phi^t\},\gamma_0,\gamma_1}(y)$
is the time needed for the flow $\{\Phi^t\}$ to
drift from the point $\gamma_0(y)$ to $\gamma_1(y)$.
That is,
$$
\Phi^{t}(\gamma_0(y))=\gamma_1(y) \mbox{ for }t=
f_{\{\Phi^t\},\gamma_0,\gamma_1}(y).$$
Then the function $f_{\{\Phi^t\},\gamma_0,\gamma_1}$ belongs to the
space 
$$E=\{f:(0,\infty)\to\R\mid f\ \ \mbox{is continuous and}
\ \ \lim_{y\to 0}f(y)=\infty\}.$$
Of course $f_{\{\Phi^t\},\gamma_0,\gamma_1}$ depends upon the
choice of $\gamma_i$. There are two umbiguities,
one coming from the parametrization of $\gamma_i$, and the other
coming from the positions of $\gamma_i$.
Let $H$ be the space of homemorphisms of $[0,\infty)$ and $C$
the space of continuous functions on $[0,\infty)$.
Define an equivalence relation $\sim$ on $E$ by
$$
f\sim f'\Longleftrightarrow f'=f\circ h+k,\ \  \exists h\in H,\
\ \exists k\in C.$$
Then clearly the equivalence class $[f_{\{\Phi^t\},\gamma_0,\gamma_1}]$
does not depend on the choice of $\gamma_i$. 
Moreover
it is an invariant of the topological conjugacy classes of
$\RR$-flows. Therefore we get a well defined map
$$\iota:\EE\to E/\sim.$$
It is easy to see that $\iota$ is injective. 
The main result of \cite{L} states that $\iota$ is surjective
as well, i.\ e.\ any $f\in E$ is realized
as $f=f_{\{\Phi^t\},\gamma_0,\gamma_1}$ for some $\RR$-flow $\{\Phi^t\}$
and curves $\gamma_i$.

The equivalence class $[f]$ of $f\in E$ is determined by how 
$f(y)$ oscilates while it tends to $\infty$ as $y\to 0$.
For example any monotone function of $E$ belongs to a single
equivalence class, which corresponds to a standard Reeb 
flow. By definition a {\em standard Reeb flow} is a flow which
is topologically conjugate to the flow $\{\varphi^t_P\}$
on $P$ given by
$$
\varphi_P^t(x,y)=(2^tx,2^{-t}y).$$

To measure the degree of oscilation of
$f\in E$, define a nonnegative valued continuous function
$f^*$ defined on $(0,1]$ by
$$
f^*(y)=\max(f\vert_{[y,1]})-f(y).
$$
Then we have the following easy lemma.

\begin{lemma} \label{easy}

(1) If $h\in H$, then $(f\circ h)^*=f^*\circ h$ in a neighbourhoof of $0$.

(2) If $k\in C$ and $y\to 0$, then $(f+k)^*(y)-f^*(y)\to 0$.

(3) There is a sequence $\{y_n\}$ converging to $0$
such that  $f^*(y_n)=0.$
\qed
\end{lemma}

For $f$ as above, define an invariant $\sigma(f)\in[0,\infty]$ by
$$\sigma(f)=\limsup_{y\to 0}f^*(y).$$
Lemma \ref{easy} implies that $\sigma(f)$ is an invariant of the class
$[f]$. We also have $\sigma(f)=0$ if and only if the class
$[f]$ is represented by a monotone function, that is, $[f]$
corresponds to a standard Reeb flow.

\section{Proof of Theorem \ref{main}}

Fix once and for all a nonstandard Reeb homeomorphism $h$ of $B$ and assume that
$h$ is the time one of a flow $\{\Phi^t\}$.
Then it can be shown that the orbit foliation $\RR$
of $\{\Phi^t\}$ is a bundle foliation in $\Int(B)$,
and therefore $\RR$ is a Reeb foliation. Let $\gamma_i$ ($i=0,1$)
and $f=f_{\{\Phi^t\},\gamma_0,\gamma_1}$ be as in Sect.\ 2.
Notice that $\sigma(f)>0$ since $\{\Phi^t\}$ must be nonstandard.

Our plan is to define ``coordinates'' of $B$ via these data,
and study the behaviour of any other flow $\{\varphi^t\}$
whose time one is $h$ using these coordinates. But it is more
convenient to work with the quotient space by $h$.
So let
$$A=B/\langle h\rangle.$$
$A$ is a non Hausdorff 2-manifold with two boundary cirles,
$\partial_iA=\partial_iB/\langle h\rangle$ ($i=0,1$).
Any neighbourhood of any point of $\partial_0A$ intersects
any neighbourhood of any point of $\partial_1A$.
Denote 
$$\Int(A)=\Int(B)/\langle h\rangle, \ \  A_i=\Int(A)\cup
\partial_iA, \ \ i=0,1.
$$
$A_i$ is a Hausdorff space homeomorphic to $S^1\times[0,\infty)$.

The flow $\{\Phi^t\}$, as well as the other flow $\{\varphi^t\}$,
induces an $S^1$-action on $A$, still denoted by the same letter.
The curve $\gamma_i$ induces a curve in $A_i$, denoted by the same
letter. One can use the parameter of the curve $\gamma_i$ as a
hight function $p$ on $A$.
Recall that by the convention of Sect.2, the points $\gamma_0(y)$ and
$\gamma_1(y)$ lie on the same leaf of $\RR$ if $y>0$.
Let us define a projection $p:A\to[0,\infty)$ as follows.
For any $\xi\in A$,
$$
p(\xi)=y\Longleftrightarrow \xi\mbox{ lies on the leaf of }\RR
\mbox{ passing through }\gamma_0(y)\mbox{ or }\gamma_1(y). $$
Of course $p(\partial_iA)=\{0\}.$ The orbit foliation of the
$S^1$ action $\{\Phi^t\}$ is now {\em horizontal}.
On the other hand we do not know what the orbit foliation
of $\{\varphi^t\}$ looks like.

Next define a
projection $\pi_i:A_i\to S^1$ as follows.
For any point $\xi\in A_i$, 
$$\pi_i(\xi)=x\Longleftrightarrow \xi=\Phi^x(\gamma_ip(\xi)).$$
Since
$$\xi=\Phi^{\pi_0(\xi)}\gamma_0p(\xi)
\ \mbox{ and }\
\xi=\Phi^{\pi_1(\xi)}\gamma_1p(\xi)
=\Phi^{\pi_1(\xi)}\Phi^{fp(\xi)}\gamma_0p(\xi),
$$
we have
\begin{equation} \label{important}
\pi_0-\pi_1=f\circ p \mbox{ mod }\Z \ \ \mbox{ \ \ on }\Int(A).
\end{equation}

There is a homeomorphism
$$\pi_0\times p:A_0\to S^1\times[0,\infty).
$$
If we use the arguments $x=\pi_0(\xi)$ and $y=p(\xi)$ on $A_0$,
 the foliation by $\pi_0$ is {\em vertical} i.\ e.\ given
by the curves $x={\rm const}.$, while
the foliation by $\pi_1$, defined on $\Int(A)$, is given by the curves
$$
x=f(y)+\mbox{const.}\mbox{ mod }\Z.$$
Both foliations are invariant by the horizontal rotation
$\{\Phi^t\}$.

Now choose a decreasing sequence of positive numbers
$$
1>y_1>y_1'>y_2>y_2'>\cdots$$
such that
\begin{equation} \label{f^*}
f^*(y_n')>\sigma(f)/2,\ \ y_n=\min\{y>y_n'\mid f^*(y)=0\},
\ \ y_n\to0\ (n\to\infty).
\end{equation}
By Lemma \ref{easy} (3), the value of the function $f^*$ is oscilating between
$0$ and around $\sigma(f)$. Therefore it is possible to choose such a
sequence.
Returning to $f$, (\ref{f^*}) implies 
\begin{equation} \label{f}
f(y_n)=\max(f\vert_{[y_n',1]}), \ \ f(y_n')<f(y_n)-\sigma(f)/2.
\end{equation}
Thus the foliation by $\pi_1$ has oscilating leaves. It is 
(topologically) tangent
to the vertical foliation at the level set $p^{-1}(y_n)$.

To prove Theorem \ref{main}, it suffices to show the existence
of a homeomorphism of $\partial_0A$ to $\partial_1A$ that conjugates
$\varphi^t\vert_{\partial_0A}$ to $\varphi^t\vert_{\partial_1A}$ for
any free $S^1$ action $\{\varphi^t\}$ on $A$. Define an $S^1$ action
$\{\varphi^t_i\}$ on $S^1$ ($i=0,1$) as the conjugate of
$\varphi^t\vert_{\partial_iA}$ by $\pi_i$, i.\ e.\ so as to satisfy
$$
\varphi_i^t\circ\pi_i=\pi_i\circ\varphi^t \mbox{ on } \partial_iA.$$
Then our goal is to show that $\varphi_0^t$ is conjugate to $\varphi_1^t$
by a homeomorphism $g$ of $S^1$ which can be chosen independently
of the $S^1$ action $\{\varphi^t\}$. But since $\{\Phi^t\}$ is one such
$S^1$ action and $\Phi_i^t$ is just a rotation,
the homeomorphism $g$ must be a rotation itself.

Besides (\ref{f}), we may assume
\begin{equation}\label{conv}
f(y_n)\to\alpha\ \mbox{ mod }\ \Z.
\end{equation}
Denote by $R_\alpha:S^1\to S^1$ the rotation by $\alpha$.
 Now our goal is to show the following proposition.

\begin{proposition} \label{beta}
For any free $S^1$-action
$\{\varphi^t\}$ on $A$, we have
$$
R_\alpha\circ\varphi_1^t=\varphi_0^t\circ R_\alpha,\ \ \forall t.$$
\end{proposition}

This follows from the following lemma.

\begin{lemma} \label{alpha}
For any free $S^1$ action $\{\varphi^t\}$ on $A$
and for any nonnegative integer $k$
$$
R_\alpha\circ\varphi_1^{1/2^k}=\varphi_0^{1/2^k}\circ R_\alpha.$$
\end{lemma}

\bigskip
The rest of this section is devoted to the proof of Lemma \ref{alpha}.
We shall first prove it for $k=1$. 
 Define a function
$\delta:\Int(A)\to\R$ by
$$
\delta(\xi)=(\pi_0\varphi^{1/2}(\xi)-\pi_1\varphi^{1/2}(\xi))
-(\pi_0(\xi)-\pi_1(\xi)).$$ 
We shall study the function $\delta$ on the circle $p^{-1}(y_n)$.
For $\xi\in p^{-1}(y_n)$, we have by (\ref{important})
\begin{equation} \label{new}
\delta(\xi)=fp\varphi^{1/2}(\xi)-f(y_n).
\end{equation}
The position of $\varphi^{1/2}(\xi)$ for $\delta(\xi)>0$
is indicated by the dot in the figure. Notice that
it must be below $p^{-1}(y_n')$.

\begin{figure}

\setlength{\unitlength}{0.7mm}
\begin{picture}(150,100)(0,0)
\put(10,80){\line(1,0){120}}
\put(100,80){\line(-4,1){70}}
\put(100,80){\line(-4,-1){100}}
\put(0,55){\line(5,-1){150}}
\put(0, 62){\line(1,0){130}}
\put(133, 80){$p^{-1}(y_n)$}
\put(133, 62){$p^{-1}(y_n')$}
\put(120,29.3){$\bullet$}
\put(0,75){\framebox(155,10)[bl]{$V_n$}}
\end{picture}
\end{figure}

There exists a {\em horizontally going point} $\xi_n(0)$ 
in $\pi^{-1}(y_n)$, i.\ e.\ a point such that
$$
p\varphi^{1/2}(\xi_n(0))=p(\xi_n(0))=y_n.
$$
For, otherwise $\varphi^{1/2}$ will
displace the curve $p^{-1}(y_n)$, sending it, say below itself.
But then $\varphi^{1/2}\circ \varphi^{1/2}$ cannot be the identity.

Notice that $\varphi^{1/2}(\xi_n(0))$ is also a horizontally going point.
By (\ref{new}), we have
\begin{equation} \label{delta}
\delta(\xi_n(0))=0\ \mbox{ and }\ \delta\varphi^{1/2}(\xi_n(0))=0.
\end{equation}
Passing to a subsequence if necessary, we may assume that
$$
\pi_i(\xi_n(0))\to \alpha_i,\ \ i=0,1.
$$
Of course we have by (\ref{important})
$$
\alpha=\alpha_0-\alpha_1.$$
For any $x\in S^1$, define $\xi_n(x)$ to be the point
on $p^{-1}(y_n)$ such that 
$$\pi_i(\xi_n(x))=\pi_i(\xi_n(0))+x\ \
(i=0,1).$$
 Then we have
\begin{equation} \label{x}
\pi_i(\xi_n(x))\to x+\alpha_i,\ \mbox{ and }\ 
\pi_i\varphi^{1/2}(\xi_n(x))\to\varphi^{1/2}_i(x+\alpha_i).
\end{equation}
To see the second assertion, notice that by the first assertion and
the fact that $p(\xi_n(x))\to 0$, 
the point $\xi_n(x)$ converges to a point $\xi_{\infty,i}(x)\in\partial_iA$.
Of cource 
$$\pi_i(\xi_{\infty,i}(x))=x+\alpha_i\ \mbox{ and } \ 
\varphi^{1/2}(\xi_n(x))\to\varphi^{1/2}(\xi_{\infty,i}(x)).
$$
Therefore
$$
\pi_i\varphi^{1/2}(\xi_n(x))\to\pi_i\varphi^{1/2}(\xi_{\infty,i}(x))
= \varphi_i^{1/2}\pi_i(\xi_{\infty,i}(x))=\varphi_i^{1/2}(x+\alpha_i),
$$
as is asserted.

Now we have
\begin{equation} \label{J} 
\delta(\xi_n(x))=(\pi_0\varphi^{1/2}(\xi_n(x))-\pi_0(\xi_n(x)))
-(\pi_1\varphi^{1/2}(\xi_n(x))-\pi_1(\xi_n(x)))
\end{equation}
$$
\longrightarrow(\varphi_0^{1/2}(x+\alpha_0)-\alpha_0)- 
(\varphi_1^{1/2}(x+\alpha_1)-\alpha_1)=J_0(x)-J_1(x),
$$
where
$$J_i=R_{\alpha_i}^{-1}\circ\varphi_i^{1/2}\circ R_{\alpha_i},$$
an involution on $S^1$.
Now (\ref{delta}) and (\ref{J}) implies that
$$
J_1(0)=J_0(0).$$
By some abuse we denote by $<$ the positive
circular order for two nearby points of $S^1$.

All we are about is to show that $J_0=J_1$. Assume for contradiction
that this is not the case.
Since $J_i$ is an involution, there is a point $0<x_0<J_i(0)$ 
such that $J_0(x_0)\neq J_1(x_0)$.
Choosing the point $x_0$ as near $0$ as we wish, we can assume
\begin{equation} \label{1/4}
\abs{J_1(x)-J_0(x)}<\sigma(f)/4\ \ \mbox{ if }\
0\leq x\leq x_0.
\end{equation}

There are two cases, one $J_0(x_0)>J_1(x_0)$
and the other $J_0(x_0)<J_1(x_0)$. But 
the latter case can be
reduced to the former case by replacing $0$ by $J_0(0)=J_1(0)$
and $x_0$ by $J_0(x_0)$.
Notice that the image by $\pi_i$ of the horizontally
going points $\varphi^{1/2}(\xi_n(0))$ converge to 
$J_0(0)=J_1(0)$. This is all we need in the argument that follows,
and therefore we can replace $0$ by $J_0(0)=J_1(0)$.

So we assume
\begin{equation} \label{case1}
J_0(x_0)> J_1(x_0).
\end{equation}

By (\ref{J}) (\ref{1/4}) and (\ref{case1}), we have for any large $n$,
\begin{equation} \label{1/3}
\abs{\delta(\xi_n(x))}<\sigma(f)/3 \mbox{ if } 0\leq x\leq x_0\ \mbox{ and }\
\delta(\xi_n(x_0))>0.
\end{equation}
Let
$$W_n=\{\xi\in A\mid fp(\xi)>f(y_n)-\sigma(f)/3\}$$
and let $V_n$ be the connected component of $W_n$ that
contains $p^{-1}(y_n)$. The subset $V_n$ is a horizontal open annulus
disjoint from $p^{-1}(y_n')$. See the figure.
By (\ref{f}), 
\begin{equation}\label{V_n}
\xi\in V_n\Longrightarrow fp(\xi)\leq f(y_n).
\end{equation}

Now we have by (\ref{new}) and (\ref{1/3})
$$
fp\varphi^{1/2}(\xi_n(x))-f(y_n)
=\delta(\xi_n(x))>-\sigma(f)/3
$$
if $0\leq x\leq x_0$. This shows that $\varphi^{1/2}(\xi_n(x))$
is contained in $W_n$. But since $\xi_n(0)$
is horizontally going, $\varphi^{1/2}(\xi_n(0))$ lies
in $V_n$. Moreover the assignment 
$$x\mapsto\varphi^{1/2}(\xi_n(x))$$
is continuous.
Thus $\varphi^{1/2}(\xi_n(x))$ lies in $V_n$ for any $0\leq x\leq x_0$.
In particular
$$
\varphi^{1/2}(\xi_n(x_0))\in V_n\ \mbox{ and }\
fp\varphi^{1/2}(\xi_n(x_0))\leq f(y_n).
$$
But the assumption $\delta(\xi_n(x_0))>0$ of (\ref{1/3}) implies
$$
fp\varphi^{1/2}(\xi_n(x_0))>f(y_n).$$
The contradiction shows that $J_0=J_1$, i.\ e.\
$$
\varphi_1^{1/2}\circ R_{\alpha}=R_{\alpha}\circ\varphi_0^{1/2}$$
for $\alpha=\alpha_0-\alpha_1$, as is required.

Now $J_0=J_1$ implies that for any large $n$
$$\xi\in p^{-1}(y_n)\Longrightarrow \abs{\delta(\xi)}<\sigma(f)/3,$$
that is, any point in $p^{-1}(y_n)$
is {\em nearly horizontally going}, meaning that it is
mapped by $\varphi^{1/2}$ into $V_n$.

\bigskip
To show Lemma \ref{alpha} for $k=2$, 
first choose a horizontally going point $\xi_n'(0)\in
p^{-1}(y_n)$
for $\varphi^{1/4}$. Its image $\varphi^{1/4}(\xi_n'(0))$ is not
horizontally going, but nearly horizontally going for $\varphi^{1/4}$.
Passing to a subsequence, we may assume
$$
\pi_i(\xi_n'(0))\to \alpha_i'.
$$
Clearly we have
$$
\alpha=\alpha_0'-\alpha_1'.$$
The point $\xi'_n(x)$ in $p^{-1}(y_n)$ is defined
just as before by 
$$\pi_i(\xi'_n(x))=\pi_i(\xi_n'(0))+x.$$

Define a function
$\delta':\Int(A)\to\R$ by
$$
\delta'(\xi)=(\pi_0\varphi^{1/4}(\xi)-\pi_1\varphi^{1/4}(\xi))
-(\pi_0(\xi)-\pi_1(\xi)),$$ 
and let
$$
J_i'=R_{\alpha_i'}^{-1}\circ\varphi^{1/4}\circ R_{\alpha_i'}.$$

Then we have
$$
\delta(\xi'_n(x))\to J_0'(x)-J_1'(x).$$

By the previous step we have shown
$$
\varphi_1^{1/2}=R_\alpha^{-1}\circ\varphi_0^{1/2}\circ R_\alpha,$$
which implies
$(J_0')^2=(J_1')^2$.
This enables us to find a point $x_0'$ playing the same role as
$x_0$ in the previous argument such that $J_0'(x_0')>J_1'(x_0')$
either near $0$ or near $J_0'(0)=J_1'(0)$.
In the latter case the point $\varphi^{1/4}(\xi_n'(0))$ is
only nearly horizontally going 
but this is enough for our purpose.
By the same argument as before, we can show $J_0'=J_1'$.

The proof for general $k$ is by an induction.

\section{Examples}

We shall construct two examples of $f\in E$.
We consider the correspoding flow $\{\Phi^t\}$ and construct the
non Hausdorff space $A$
as in Sect.3. Properties of examples are stated in terms of
the $S^1$ action on $A$.
All the notations of Sect.\ 3
will be used.

\begin{example} \label{e1}
There exists $f\in E$ such that $\sigma(f)=1$ satisfying the
following property: 
For any $S^1$ action $\{\psi^t\}$ on $S^1$, there
is an $S^1$ acion $\{\varphi^t\}$ on $A$ such that
$\varphi_i^t=\psi^t$ ($i=0,1$).
\end{example}

The construction of $f$ goes as follows.
Let
$$
y_1>y_1'>y_2>y_2'\cdots$$
be a sequence converging to 0. Define $f$ such that
$$
f(y_n)=n,\ \mbox{ and }\ f(y_n')=n-1$$
and that $f$ is monotone on the complementary
intervals.

On the circles $p^{-1}(y_n)$ and $p^{-1}(y_n')$,
$\pi_0=\pi_1$ mod $\Z$. The desired flow is to preserve these
circles and to be the conjugate of $\psi^t$
by $\pi_i$ there. The complementary regions are open annulus, and
there the foliations by $\pi_0$ and $\pi_1$ are transverse,
thanks to the monotonicity assumption on $f$. Therefore
one can define $\varphi^t$ so as to satisfy
$
\pi_i\circ\varphi^t=\psi^t\circ\pi_i$ ($i=0,1$).

\begin{example}
There exists $f\in E$ such that any $S^1$ action $\{\varphi^t\}$ on $A$
satisfies $\varphi_i^t=R_t$ ($i=0,1$).
\end{example}

Using the sequence of Example \ref{e1}, define $f\in E$ such that
$$
f(y_n)=n\beta\ \mbox{ and }\ f(y_n')=n\beta-1,$$
for some irrational $\beta>0$
and that $f$ is monotone on the complementary intervals.
Then
$$
\pi_0-\pi_1=n\beta\ \mbox{ on }\ p^{-1}(y_n).$$
Let $\{\varphi^t\}$ be an arbitrary $S^1$ action on $A$.
Then any point $\tau\in S^1$ is an accumulation point of
$\pi_0(\xi_n(0))-\pi_1(\xi_n(0))$. The argument of Sect.\ 3
shows that 
$$\varphi_1^t\circ R_{\tau}=R_{\tau}\circ \varphi_0^t,\ \ \forall t,\tau.
$$
This clearly shows that 
$$
\varphi_0^t=\varphi_1^t=R_t.$$

\end{document}